\documentclass[12pt,draft]{article}

\usepackage[english]{babel}

\usepackage{amssymb}
\usepackage{latexsym}
\usepackage[round]{natbib}
\usepackage{url}

\usepackage[final]{graphicx}

\newcommand{\prob}{\mathrm{Pr}}

\newcommand\cvd{\hfill\raisebox{-0.5ex}{$\Box$}}
\newcommand\cvdbis{\hfill\raisebox{1.5ex}{$\Box$}}

\newcommand{\bc}{\begin{center}}
\newcommand{\ec}{\end{center}}

\newcommand{\bit}{\begin{itemize}}
\newcommand{\eit}{\end{itemize}}

\newcommand{\be}{\begin{eqnarray*}}
\newcommand{\ee}{\end{eqnarray*}}
\newcommand{\ben}{\begin{eqnarray}}
\newcommand{\een}{\end{eqnarray}}

\newcommand{\ntr}{n_{\triangle}}
\newcommand{\xtr}{x_{\triangle}}
\newcommand{\atr}{a_{\triangle}}

\newcommand{\nplus}{n_{++}}
\newcommand{\xplus}{x_{++}}
\newcommand{\aplus}{a_{++}}

\newcommand{\noi}{n_{01}}
\newcommand{\xoi}{x_{01}}
\newcommand{\aoi}{a_{01}}

\newcommand{\nio}{n_{10}}
\newcommand{\xio}{x_{10}}

\newcommand{\noo}{n_{00}}

\usepackage{layout}

\newtheorem{definition}{Definition}[section]

\newtheorem{lemma}{Lemma}[section]
\newtheorem{prop}{Proposition}[section]

\newtheorem{example}{Example}[section]

\begin{document}

\title{Intrinsic tests for the equality\\
of two correlated proportions}
\author{
Guido Consonni\\
University of Pavia\\
\url{guido.consonni@unipv.it}
\and
Luca La Rocca\\
University of Modena and Reggio Emilia\\
\url{luca.larocca@unimore.it}
}

\maketitle


\begin{abstract}
Correlated proportions arise in longitudinal (panel) studies.
A typical example is the ``opinion swing'' problem:
\emph{``Has the proportion of people favoring a politician changed
after his recent speech to the nation on TV?''}.
Since the same group of individuals is interviewed before and after the speech,
the two proportions are correlated.
A natural null hypothesis to be tested is
whether the corresponding population proportions are equal.
A standard Bayesian approach to this problem
has already been considered in the literature,
based on a Dirichlet prior for the cell-probabilities
of the underlying two-by-two table under the alternative hypothesis,
together with an induced prior under the null.
In lack of specific prior information, a diffuse (e.g.~uniform) distribution may be used.
We~claim that this approach is not satisfactory, since in a testing problem
one should make sure that the prior under the alternative be adequately centered
around the region specified by the null,
in order to obtain a fair comparison between the two hypotheses.
Following an intrinsic prior methodology, we develop two strategies
for the construction of a collection of objective priors
increasingly peaked around the null.
We provide a simple interpretation of their structure
in terms of weighted imaginary sample scenarios.
We illustrate our method by means of three examples,
carrying out sensitivity analysis
and providing comparison with existing results.
\bigskip

\noindent\textit{Keywords}:
Bayes factor; Bayesian robustness, Dirichlet prior;
marginal homogeneity; matched pair;
objective Bayes.
\end{abstract}

\section{Introduction}

Panel data are often used to investigate changes in opinion or behavior.
For instance, suppose the same group of individuals is interviewed about
their support for the President before and after the State of the Union Address.
A natural issue is whether the speech has generated or not a net change in attitude.
If the answer can only be ``I support the President''
or ``I do not support the President'', i.e.~it is binary,
the problem reduces to the comparison of two proportions in the population:
of those supporting the President  before  and after the address.
We say that no net change has occurred if the two proportions are equal.
Clearly, the two sample proportions are correlated,
since the same group of people is involved in the two measurements.
A Multinomial model involving the four pairs of outcomes
``(Support, Support)'', ``(Support, Not Support)'',
``(Not Support,  Support)'', ``(Not Support, Not Support)''
is the appropriate sampling scheme to be used.

A classic frequentist test dates back to \citet{McNe:1947}.
The first Bayesian analysis is attributed to \citet{Alth:1971},
while more recent contributions include \citet{Broe:Greg:1996},
\citet{Iron:Pere:Tiwa:2000} and \citet{Kate:Papa:Dell:2001}.
In particular, the latter two papers compute the Bayes Factor (BF)
for the hypothesis of no net change \textit{versus} the alternative
that a net change has occurred, starting from a ``default'' Dirichlet prior
on the unconstrained cell-probabilities.
In~lack of specific prior information, a diffuse (e.g.~uniform) distribution
may be used. The prior distribution under the hypothesis of no net change
is then induced from the Dirichlet prior.
We call this a ``default''~BF for testing the equality of two correlated proportions,
and claim that it is not satisfactory, since it might unduly favor the null hypothesis.
To overcome this drawback, we~keep the same prior under the null
and~construct a suitable intrinsic prior \citep{Berg:Peri:1996} under the alternative,
using an objective approach along the lines recently described
by \cite{Case:More:2006,Case:More:2007}.
The intrinsic prior is centered around the region where the null obtains,
and the corresponding BF results in a fairer comparison between the two hypotheses.
As the amount of centering is controlled by a scalar hyper-parameter,
we actually have a family of prior distributions.
We suggest studying BF sensitivity to prior choice within this family.

The rest of the paper is organized as follows.
In Section~2 we give a new derivation of the above default BF for testing the
equality of two correlated proportions;
in Section~3 we discuss the concept of an intrinsic prior for testing hypotheses,
and derive intrinsic priors for the problem at hand;
in~Section~4 we illustrate our method by means of three examples.
Finally, Section~5 contains a brief discussion.

\section{Default approach}

We assume the data consist of $n_{++}$ observations on a pair
of binary  variables, $(U,V)$ say, each taking values in $\{ 0,1
\}$, with $n_{++}$ known by design. To fix ideas,
let  $U=1$ ($V=1$) represent the opinion
``I support the President'', and~$U=0$ ($V=0$) the opinion
``I do not support the President'', before (after) the State of the Union Address.
A  swing occurs whenever $(U=0,V=1)$ or $(U=1,V=0)$.
We can arrange the four counts in a $2 \times 2$ table
$n=(n_{00},n_{01},n_{10},n_{11})$,
where $n_{00}+n_{01}+n_{10}+n_{11}=n_{++}$.
If the data are independent, conditionally on some common parameter,
the statistical model for the data is quadrinomial
with matrix of cell-probabilities $\pi=(\pi_{00},\pi_{01},\pi_{10},\pi_{11})$,
where~$\pi_{00}+\pi_{01}+\pi_{10}+\pi_{11}=1$.
If no net change has occurred, the \textit{marginal} distributions of~$U$ and~$V$
are the same, so~that the problem under investigation consists in
testing \textit{marginal homogeneity} in a two-way contingency table,
i.e.~$\pi_{0+}=\pi_{+0}$, where subscript ``+'' denotes summation
over the replaced suffix. In a $2 \times 2$ table, however,
marginal homogeneity reduces to \textit{symmetry}: $\pi_{01}=\pi_{10}$.

We take advantage of the following reparametrization:
$\eta = \pi_{01}+\pi_{10}$, $\theta = \frac{\pi_{01}}{\pi_{01}+\pi_{10}}$
and $\gamma = \frac{\pi_{00}}{\pi_{00}+\pi_{11}}$.
Notice that $\eta$ is the (unconditional) probability of observing a swing,
whereas $\theta$ is the conditional probability of a swing from~0 to~1,
given that a swing occurs,
and $\gamma$ is the conditional probability of no support for the President,
given that no swing occurs. The sampling distribution of~$n$ can thus be written as
\begin{eqnarray}
f(n|\eta,\theta,\gamma)
& = &
g(\ntr|\eta)h(n_{01}|\ntr,\theta)k(n_{00}|\ntr,\gamma),
\label{eq:quadrinomial}
\end{eqnarray}
writing $\ntr=n_{01}+n_{10}$ for the sum of the off-diagonal counts and letting
$g(\ntr|\eta)=\mbox{Bin}(\ntr|\nplus,\eta)$,
$h(\noi|\ntr,\theta)=\mbox{Bin}(\noi|\ntr,\theta)$ and
$k(\noo|\ntr,\gamma)=\mbox{Bin}(\noo|\nplus-\ntr,\gamma)$,
where $\mbox{Bin}(y|m,\psi)$ is
the Binomial probability function with $m$ trials,
and probability of success~$\psi$, evalutated at~$y$.

With reference to the $(\eta,\theta,\gamma)$-parametrization,
the null hypothesis of no net change can be formulated as
$H_0:\theta=\frac{1}{2}$, whereas the alternative hypothesis is
$H: \theta \neq \frac{1}{2}$.
Notice that $\eta$ and $\gamma$ are nuisance parameters,
and that the parameter of interest $\theta$ only appears
in $h(n_{01}|n_{\triangle},\theta)$,
the conditional distribution of~$\noi$ given~$\ntr$.
Accordingly, \citet{Iron:Pere:Tiwa:2000} base their analysis exclusively on
the partial sampling model $h(n_{01}|n_{\triangle},\theta)$.
This is also true for the classical test by \cite{McNe:1947},
which only makes use of the off-diagonal individual swing-counts $(n_{01},n_{10})$,
and for \citet{Alth:1971}, who performs Bayesian inference
on the difference $\pi_{01}-\pi_{10}$.
However, this issue is controversial, as some Authors suggest that
the whole sample size $n_{++}$ should also be taken into
consideration; see for instance \citet[p.~442]{Agre:2002}.
Therefore, in the following,
we consider the full sampling model~(\ref{eq:quadrinomial}).
Eventually, neither $g(\ntr|\eta)$ nor $k(n_{00}|\ntr,\gamma)$
will play a role in the default approch, but $g(\ntr|\eta)$
will do in the intrinsic approach. We shall return to this issue
later on in the paper, when considering the specification of
an intrinsic prior, and in the Discussion.

As in previous Bayesian analyses of this problem,
we assign a Dirichlet prior on the matrix of cell-probabilities~$\pi$.
We do this for convenience, and to ease comparisons with existing results.
In this way, we shall obtain a default BF
for testing the equality of two correlated proportions.
Later on, however, we shall significantly depart from  previous works
on this topic, because the prior under the alternative hypothesis will be
suitably adjusted to the testing problem under consideration.
Next lemma shows that the resulting prior on~$(\eta,\theta,\gamma)$
is a product of independent beta distributions,
shedding light on the implications of choosing a Dirichlet prior for~$\pi$.

\begin{lemma}
\label{lem:prior}
Let $\pi \sim \mbox{Dir}(a)$, where $a=(a_{00},a_{01},a_{10},a_{11})$
is a two-by-two matrix of (strictly) positive values.
Then:
i) $\eta \sim \mbox{Beta}(a_{01}+a_{10},a_{00}+a_{11} )$;
ii)  $\theta \sim \mbox{Beta}(a_{01},a_{10})$;
iii) $\gamma \sim \mbox{Beta}(a_{00},a_{11})$;
iv) $\eta$, $\theta$ and $\gamma$ are mutually independent.
\end{lemma}
\textit{Proof}.
Substitute $\pi_{00}=\gamma(1-\eta)$, $\pi_{01}=\theta\eta$
and $\pi_{10}=(1-\theta)\eta$ in the Dirichlet density of~$(\pi_{00},\pi_{01},\pi_{10})$,
and multiply by Jacobian $\eta(1-\eta)$ to get the density of $(\eta,\theta,\gamma)$.
\cvd
\smallskip


Recall that the BF of $H$ \emph{versus} $H_0$, or more briefly in favor of~$H$,
having observed~$n$, is~given by $BF_{H,H_0}(n)=\frac{m_H(n)}{m_{H_0}(n)}$,
where $m_H(n)$ and $m_{H_0}(n)$ are the marginal probabilites of~$n$
under $H$ and $H_0$, respectively. More specifically, we have
$
m_H(n) =
\int\!\!\!\int\!\!\!\int
f(n|\eta,\theta,\gamma)p_H(\eta,\theta,\gamma)\,
d\eta d\theta d\gamma
$
and
$
m_{H_0}(n) =
\int\!\!\!\int
f(n|\eta,\theta=1/2,\gamma)p_{H_0}(\eta,\gamma)\,
d\eta d\gamma
$,
where $p_H(\eta,\theta,\gamma)$ is given by Lemma~\ref{lem:prior}
and $p_{H_0}(\eta,\gamma)$ has yet to be specified.
If we derive it from $p_H(\eta,\theta,\gamma)$ as
i) the marginal $p_H(\eta,\gamma)$
or ii) the conditional $p_H(\eta,\gamma|\theta=1/2)$
we obtain the same result:
a product of two independent beta densities.
Using this induced prior on~$H_0$,
the following lemma produces a default BF \emph{against}
the equality of two correlated proportions.

\begin{lemma}
\label{lem:defaultBF}
If $\pi \sim \mbox{Dir}(a)$, then  a default BF in favor of~$H$ is given by
\begin{eqnarray}
BF_{H,H_0}(n)=\frac{2^{n_{\triangle}}B(a_{01}+n_{01},a_{10}+n_{10})}{B(a_{01},a_{10})},
\label{eq:defaultBF}
\end{eqnarray}
where $B(a,b)=\int_{0}^{1}t^{a-1}(1-t)^{b-1}dt$, $a>0$, $b>0$,
is the Beta special function.
\end{lemma}
\textit{Proof}.
Since the prior density of $(\eta,\gamma)$ under~$H_0$
corresponds to~$p_H(\eta,\gamma|\theta=1/2)$, the Savage density ratio applies;
see e.g.~\citet[Sec.~7.16]{Ohag:Fors:2004}. This gives
$BF_{H,H_0}(n)=\frac{p_H(\theta=1/2)}{p_H(\theta=1/2|n)}$,
where $p_H(\theta|n)$ is the posterior density of~$\theta$.
Now, from ii) of Lemma \ref{lem:prior},
the prior  of $\theta$ is $\mbox{Beta}(a_{01},a_{10})$, whence
$p_H(\theta=1/2)=\frac{1}{B(a_{01},a_{10})} \left(\frac{1}{2} \right)^{a_{01}+a_{10}}$;
on the other hand, using (\ref{eq:quadrinomial}) and the independence of
$\eta$, $\theta$ and $\gamma$, the posterior density of $\theta$ is easily seen
to be $\mbox{Beta}(a_{01}+n_{01},a_{10}+n_{10})$.
The result follows.
\cvd
\smallskip

Lemma~\ref{lem:defaultBF} is also derived in~\cite{Kate:Papa:Dell:2001},
by direct computation of~$m_H(n)$ and~$m_{H_0}(n)$, and in \citet{Iron:Pere:Tiwa:2000},
using the conditional model $h(n_{01}|n_{\triangle},\theta)$~in~(\ref{eq:quadrinomial}).


\section{Intrinsic approach}

\subsection{General}
\label{subs:general}

The BF is notoriously sensitive to prior specifications.
This suggests that the choice of prior distributions should be
especially thoughtful. In~particular, we~would like to consider one specific
aspect of prior specification which is often overlooked.
When testing $H_0$ \emph{versus} $H$,
the $H_0$ hypothesis must have some reasonable grounds,
otherwise we would not have considered it in the first place.
Assuming that $H_0$ is nested in~$H$, this~suggests that the prior under $H$
should be centered around the region in which $H_0$ obtains.
Otherwise, when data are just in reasonable accord with $H_0$,
the prior under~$H$ is simply  ``wasting away'' prior probability mass
in regions that are too unlikely to be supported,
and~$H_0$ would be unduly favored.
This aspect is strictly related to the Jeffreys-Lindley's paradox;
see \citet[p.~234]{Robe:2001}.
Notice that centering the prior under $H$ around values consistent with~$H_0$
does not lend support to~$H_0$, but rather to~$H$:
a~point which is often overlooked in Bayesian applications.
The previous considerations are lucidly spelled out
in \cite{Case:More:2007}; see also  \cite{Case:More:2006}.
We now briefly review their approach, showing some basic facts about
intrinsic priors and the associated BFs.

Let $y$ be the actual data,
and denote with $x$ a corresponding set of ``imaginary'' data.
We shall assume, for simplicity of exposition, that the data are discrete,
as in the problem that motivates our paper.
Let the sampling distribution of~$y$ be $f_H(y|\lambda_H)$, under~$H$,
and $f_{H_0}(y|\lambda_{H_0})$, under $H_0$.
Let the priors under both models be given, and denote them with $p_H(\lambda_H)$
and $p_{H_0}(\lambda_{H_0})$, respectively.
We shall not insist on the nature of these priors,
although they will typically be standard ``estimation-based'' priors.
In~this paper, we shall restrict our attention to proper priors;
see \cite{More:Bert:Racu:1998} for an intrinsic limiting procedure
dealing with improper priors. Finally, let
$m_H(x)=\int f_H(x|\lambda_H)p_H(\lambda_H)d\lambda_H$
and
$m_{H_0}(x)=\int f_{H_0}(x|\lambda_{H_0})p_{H_0}(\lambda_{H_0})d\lambda_{H_0}$
be the marginal probabilities of the imaginary data $x$
under $H$ and $H_0$, respectively.

\begin{definition}
\label{def:intrinsicPrior}
The intrinsic prior for $\lambda_H$ conditionally on $H_0$ is given by
$$
p^I_{H}(\lambda_H|H_0) =
p_H(\lambda_H)E_{\lambda_H}\left[\frac{m_{H_0}(x)}{m_H(x)}\right] =
p_H(\lambda_H)E_{\lambda_H}\left[BF_{H,H_0}(x)^{-1}\right],
$$
where $BF_{H,H_0}(x)$ is the BF of $H$ \emph{versus} $H_0$
having observed the imaginary data $x$, and the expectation is taken
with respect to the sampling distribution $f_H(x|\lambda_H)$.
\end{definition}


\begin{prop}
\label{prop:integratesToOne}
If $p_H$ and $p_{H_0}$ are proper, the intrinsic prior $p_H^I$ is also proper.
\end{prop}
\textit{Proof}.
Interchanging the sum over~$x$ and the integral with respect to~$\lambda_H$, we find
$$
\int p^I_{H}(\lambda_H|H_0) d\lambda_H =
\sum_x \frac{m_{H_0}(x)}{m_H(x)}
\int f_{\lambda_H}(x|\lambda_H) p_H(\lambda_H) d\lambda_H =
1.
$$
\cvdbis
\vspace{-3.5ex}
\begin{prop}
The intrinsic prior is equivariant to reparametrization:
if $\lambda_H\!\!\mapsto\!\!\psi_H$, then
$p^{I,\psi_H}_{H}(\psi_H|H_0)=p^{I,\lambda_H}_{H}(\lambda_H(\psi_H)|H_0)
\times |J_{\lambda_H}(\psi_H)|,
$
where $J_{\lambda_H}(\psi_H)$ is the Jacobian.
\end{prop}
\textit{Proof}.
Expressing the sampling distribution of~$y$ in terms of~$\lambda_H$, we find
$$
p^{I,\psi_H}_{H}(\psi_H|H_0)=
E_{\lambda_H(\psi_H)}\left[\frac{m_{H_0}(x)}
{m_H(x)}\right]p^{\lambda_H}_{H}(\lambda_H(\psi_H))
\times |J_{\lambda_H}(\psi_H)|
$$
and the thesis follows.
\cvd
\smallskip

The two propositions above show that Definition \ref{def:intrinsicPrior}
is well-posed. Note, however, that the intrinsic prior depends on the sampling
distribution of~$x$, and in particular on its sample size (prior sample size).
As for the interpretation of Definition~\ref{def:intrinsicPrior},
the following proposition is enlightening.
\begin{prop}
\label{prop:pIasMixture}
The intrinsic prior is a mixture of imaginary posteriors:
\begin{eqnarray}
p^I_{H}(\lambda_H|H_0) & = & \sum_x  p_H(\lambda_H|x) m_{H_0}(x).
\label{eq:pIasMixture}
\end{eqnarray}
\end{prop}
\textit{Proof}.
By direct computation, we find
$$
p^I_{H}(\lambda_H|H_0)=
p_H(\lambda_H)E_{\lambda_H}\left[\frac{m_{H_0}(x)}{m_H(x)}\right]=
\sum_x \frac{m_{H_0}(x)}{m_H(x)} p_H(\lambda_H) f_{H}(x|\lambda_H)
$$
and (\ref{eq:pIasMixture}) follows.
\cvd

Proposition~\ref{prop:pIasMixture} shows that
the intrinsic prior will be more dominated
by posteriors corresponding to imaginary realizations that are more supported
under~$H_0$, and this explains why it will be more concentrated
around values consistent with $H_0$ than the original prior $p_H(\lambda_H)$.
The extent of this concentration will depend both on the structure of $p_H(\lambda_H)$
and on prior sample size. Notice that $p^I_{H}(\lambda_H|H_0)$
is an example of~\emph{expected posterior prior} (with base model $H_0$)
as defined by~\cite{Pere:Berg:2002}.
As a consequence of Proposition \ref{prop:pIasMixture}, intrinsic
prior BFs will also be mixtures of conditional BFs, over all
possible configurations of the imaginary data $x$, with weights
given by the marginal distribution of $x$ under $H_0$.

\begin{definition}
For data $y$, the intrinsic prior BF in favor of $H$ is given by
\begin{eqnarray}
BF^I_{H,H_0}(y)=\frac{m^I_H(y)}{m_{H_0}(y)},
\label{eq:BFIdef}
\end{eqnarray}
where $m^I_H(y)=\int f_H(y|\lambda_H) p^I_{H}(\lambda_H|H_0)$
is the intrinsic marginal probability of~$y$.
\end{definition}

\begin{prop}
\label{prop:BFIasMixture}
For data $y$, the intrinsic prior BF in favor of $H$ is given by
\begin{eqnarray}
BF^I_{H,H_0}(y)=\sum_x BF_{H,H_0}(y|x) m_{H_0}(x),\nonumber
\end{eqnarray}
where
$
BF_{H,H_0}(y|x)
\quad=\quad
\frac{\int f_{H}(y|\lambda_H) p_H(\lambda_H|x)d\lambda_H}{m_{H_0}(y)}
\quad=\quad
\frac{m_H(y|x)}{m_{H_0}(y)}.
$
\end{prop}
\textit{Proof}.
Substituting (\ref{eq:pIasMixture}) in~(\ref{eq:BFIdef}),
then interchanging the sum and the integral, we obtain
\begin{eqnarray*}
BF^I_{H,H_0}(y) & = &
\frac{1}{m_{H_0}(y)} \int f_H(y|\lambda_H)
\sum_x p_H(\lambda_H|x) m_{H_0}(x) d\lambda_H\\
& = &
\sum_x m_{H_0}(x)
\frac{\int f_{H}(y|\lambda_H) p_H(\lambda_H|x) d\lambda_H}{m_{H_0}(y)}
\end{eqnarray*}
and the thesis follows.
\cvd


\subsection{Conditionally-Intrinsic Prior and Test}
\label{subs:single-intrinsic}
The default BF of Lemma~\ref{lem:defaultBF}
for testing the equality of two correlated proportions
was obtained in the full sampling model~(\ref{eq:quadrinomial}),
although it can also be derived using the conditional model
$h(\noi|\ntr,\theta)$. The following argument provides
a Bayesian justification for basing inference solely
on the conditional model.

Since $\eta$ and $\gamma$ are nuisance parameters
in the testing problem under consideration, a natural suggestion is
first to integrate them out, thus obtaining an ``integrated model'' which only depends
on the parameter of interest $\theta$, then to proceed with the testing problem.
Consider first the alternative hypothesis $H$.
Using (\ref{eq:quadrinomial}) and the independence of $\eta$, $\theta$ and~$\gamma$,
we immediately obtain that the integrated model under~$H$ is given by
\begin{eqnarray}
\label{eq:integratedModelH}
f^*_H(n|\theta)=m_H(\ntr)m_H(n_{00}|\ntr)h(n_{01}|\ntr,\theta),
\end{eqnarray}
where $m_H(n_{\triangle})= \int g(n_{\triangle}|\eta)p_H(\eta)d\eta$
and $m_H(n_{00}|n_{\triangle})= \int k(n_{00}|\ntr,\gamma)p_H(\gamma)d\gamma$.
Consider next the null hypothesis $H_0$,
whose sampling model is obtained from~(\ref{eq:quadrinomial})
by setting $\theta=1/2$.
To integrate out $\eta$ and $\gamma$ we need a prior $p_{H_0}(\eta,\gamma)$.
Since they are nuisance parameters, one possible argument is that this
prior should be the same as that under~$H$;
see e.g.~\citet[Sects.~11.29-11.33]{Ohag:Fors:2004} in a different context.
In~other words $p_{H_0}(\eta,\gamma)$ should be obtained from $p_H(\eta,\theta,\gamma)$
by \textit{marginalization}.
Alternatively, one could think of obtaining $p_{H_0}(\eta,\gamma)$
by \textit{conditioning},
i.e.~$p_{H_0}(\eta,\gamma)=p_{H}(\eta,\gamma|\theta=1/2)$.
In either case we obtain the same result,
because of the independence of $\eta$, $\theta$ and $\gamma$,
so that the integrated model under~$H_0$ is
\begin{eqnarray}
f^*_{H_0}(n)=m_H(n_{\triangle})m_H(n_{00}|\ntr)
h\left( n_{01}|\theta=1/2,n_{\triangle} \right),\nonumber
\end{eqnarray}
where $m_H(n_{\triangle})$ and $m_H(n_{00}|\ntr)$ are the same as
in (\ref{eq:integratedModelH}).
From the integrated models $f^*_{H}(n|\theta)$ and $f^*_{H_0}(n)$
one can now compute the BF, by taking the ratio
$\frac{\int f^*_H(n|\theta) p_H(\theta) d\theta}{f^*_{H_0}(n)}$.
We obtain $\frac{\int h(n_{01}|n_{\triangle},\theta)p_H(\theta) d\theta}
{h(n_{01}|n_{\triangle},\theta=1/2)}$, as factors $m_H(n_{\triangle})$
and $m_H(n_{00}|\ntr)$ clearly cancel out.
This is the default BF based on the \textit{conditional} sampling model
$h(n_{01}|n_{\triangle},\theta)$.
Hence, provided we are willing to take the same marginal prior for $(\eta,\gamma)$
under $H$ and $H_0$, we can forget about $g(\ntr|\eta)$ and $k(\noo|\ntr,\gamma)$
in~(\ref{eq:quadrinomial}).

In the following of this subsection, we apply the intrinsic prior methodology
based on the conditional sampling model $h(n_{01}|n_{\triangle},\theta)$:
the resulting prior will be called ``Conditionally-Intrinsic'' (CI).
Both its density and the corresponding BF will be labeled with a superscript ``CI'';
moreover, for the sake of clarity, other intermediate quantities will have
the superscript ``Co'' to remind the reader that they all refer to the above
Conditional model. The general notation of Subsection \ref{subs:general}
will be reserved for the next subsection,
wherein the full sampling model~(\ref{eq:quadrinomial}) will be employed.
Note that, within the CI-procedure, the imaginary data $x$ are
represented by $(\xoi,\xtr)$, where  $\xtr=\xoi+\xio$,
or equivalently by $(\xoi,\xio)$.

We start by defining  the CI-prior, which is given by
\begin{eqnarray}
p^{CI}_{H}(\theta|H_0)=p_H(\theta)E_{\theta}\left[BF^{Co}_{H,H_0}(x)\right]^{-1},
\nonumber
\end{eqnarray}
where $
BF^{Co}_{H,H_0}(x)=
\frac{\int h(x_{01}|x_{\triangle},\theta) p_H(\theta) d\theta}
{h(x_{01}|x_{\triangle},\theta=1/2)}.\nonumber
$
\begin{prop}
\label{prop:CIpriorForTheta}
The CI-prior for $\theta$ conditionally on $H_0$ is given by
\begin{eqnarray}
p^{CI}_{H}(\theta|H_0)=\sum_{x_{01}=0}^{x_{\triangle}}
{x_{\triangle} \choose x_{01}}
\left (\frac{1}{2} \right )^{x_{\triangle}}
{\mbox Beta}(\theta|a_{01}+x_{01},\atr -a_{01}+x_{\triangle}-x_{01}),\nonumber
\end{eqnarray}
where $\mbox{Beta}(\theta|A,B)$ is the Beta density
with parameters $A$ and $B$, evaluated at~$\theta$.
\end{prop}

\textit{Proof}. The general structure of the intrinsic prior
is given in (\ref{eq:pIasMixture}).
It is immediate to check that $p^{Co}_{H}(\theta|x)$ is actually
${\mbox Beta}(\theta|a_{01}+x_{01},\atr-\aoi+x_{\triangle}-x_{01})$.
On the other hand, it holds that
$m^{Co}_{H_0}(x)=h(x_{01}|x_{\triangle},\theta=1/2)=
{x_{\triangle}\choose x_{01}}\left(\frac{1}{2} \right )^{x_{\triangle}}$.
\cvd
\smallskip

\begin{figure}
\begin{center}
\includegraphics[width=3.3in]{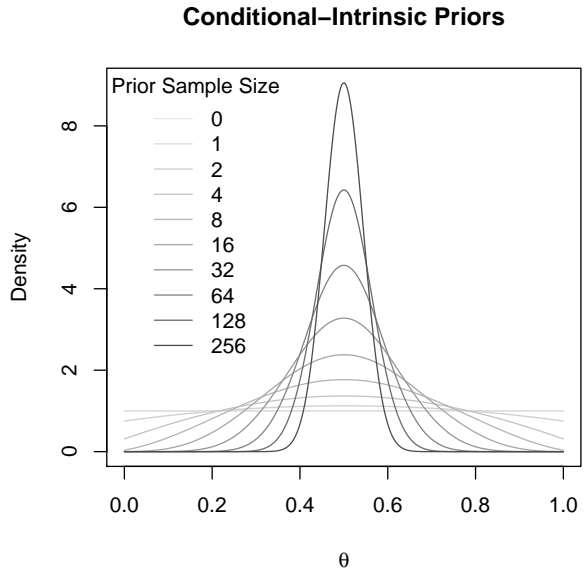}
\end{center}
\caption{Conditional-intrinsic prior densities, for ten different prior sample sizes,
starting from a uniform distribution on the cell-probabilities.}
\label{fig:CI prior}
\end{figure}

%
Proposition~\ref{prop:CIpriorForTheta} shows that the CI-prior for
$\theta$ is  a finite mixture of Beta distributions with
$x_{\triangle}+1$ components, where $x_{\triangle}$ is the prior
swing count. Figure \ref{fig:CI prior} plots the CI-prior for
$\theta$, starting from a uniform distribution on $\pi$, for ten
different prior sample sizes (values of~$\xtr$). We include the
zero sample size, which gives back the uniform default prior for
$\theta$. Notice how the CI-prior gets more and more peaked around
the null value $\theta=1/2$, as prior sample size grows.

We now turn to the derivation of the CI-prior BF, namely
\begin{eqnarray}
\label{eq:BFSI}
BF^{CI}_{H,H_0}(n)=\frac{m^{CI}_H(n)}{m^{Co}_{H_0}(n)},\nonumber
\end{eqnarray}
where $m^{Co}_{H_0}(n)=h(n_{01}|n_{\triangle},\theta=1/2)$ and
$m^{CI}_H(n)=\int h(n_{01}|n_{\triangle},\theta)
p^{CI}_H(\theta|H_0)d\theta$.

\begin{prop}
\label{prop:CIpriorBF}
The CI-prior BF in favor of~$H$ is given by
\begin{eqnarray}
BF^{CI}_{H,H_0}(n)=\sum_{x_{01}=0}^{x_{\triangle}}
{x_{\triangle} \choose x_{01}}
\left (\frac{1}{2} \right )^{x_{\triangle}}\nonumber
BF^{Co}_{H,H_0}(n|x)
\end{eqnarray}
where
$
BF^{Co}_{H,H_0}(n|x)
=
2^{n_{\triangle}}\frac{
B(a_{01}+x_{01}+n_{01},\atr-a_{01}+x_{\triangle}-x_{01}+n_{\triangle}-n_{01})}{B(a_{01}+x_{01},\atr-a_{01}+x_{\triangle}-x_{01})}.
$
\end{prop}
\textit{Proof}.
From Proposition \ref{prop:BFIasMixture} we immediately derive that the CI-prior BF
is a mixture of conditional BFs. The weights of the mixture are exactly those
of Proposition~\ref{prop:CIpriorForTheta}. On the other hand
$
BF^{Co}_{H,H_0}(n|x)
=\frac{\int h(n_{01}|n_{\triangle},\theta)
p_H^{Co}(\theta|x)d\theta}{h(n_{01}|n_{\triangle},\theta=1/2)},
$
with $p_H(\theta|x)$  a
$\mbox{Beta}(\theta;a_{01}+x_{01},\atr-a_{01}+x_{\triangle}-\xoi)$ distribution.
Standard calculations lead to the result.
\cvd
%
%

\subsection{Intrinsic Prior and Test}
In this subsection we work in the full sampling model (\ref{eq:quadrinomial}),
involving all three parameters $\eta$, $\theta$ and $\gamma$,
and derive the corresponding intrinsic prior. This is defined~as
\begin{eqnarray}
p^I_{H}(\eta,\theta,\gamma|H_0)
&=&
p_H(\eta,\theta,\gamma)E_{\eta,\theta,\gamma}\left[BF_{H,H_0}(x)^{-1}\right],\nonumber
\end{eqnarray}
according to Definition~\ref{def:intrinsicPrior},
where $BF_{H,H_0}(x)$ is given by~(\ref{eq:defaultBF}) with~$x$ in place of~$n$.
Since $\eta$, $\theta$ and $\gamma$ are independent,
and $BF_{H,H_0}(x)$ does not depend on~$x_{00}$,
it~turns out that
$p^I_{H}(\eta,\theta,\gamma|H_0)=p_H(\gamma)p^I_{H}(\eta,\theta|H_0)$,
where
$p^I_{H}(\eta,\theta|H_0)=
p_H(\eta,\theta)E_{\eta,\theta}\left[BF_{H,H_0}(x)^{-1}\right]$.
In other words, under the I-prior $\gamma$ is still independent of $(\eta,\theta)$,
and its marginal is unchanged. As a consequence, parameter~$\gamma$ and data~$x_{00}$
play no role and can be disregarded.
In the following we shall focus our attention on parameters $\eta$ and $\theta$,
and imaginary data~$x$ will be represented by the triple
$(\xoi,\xio,\xplus)$, or equivalently by~$(\xoi,\xtr,\xplus)$, where $\xplus$
is the prior sample size.

\begin{prop}
\label{prop:Iprior}
The intrinsic prior for $(\eta,\theta)$ conditionally on $H_0$ is given by
\begin{eqnarray}
p_H^{I}(\eta,\theta|H_0) & = &
\sum_{\xoi=0}^{\xplus}\sum_{\xtr=\xoi}^{\xplus}
Beta(\eta|\atr+\xtr,\aplus-\atr+\xplus-\xtr) \nonumber \\
& \times &
Beta(\theta|a_{01}+x_{01},\atr-a_{01}+\xtr-x_{01}) \nonumber \\
 & \times &
{x_{\triangle}\choose x_{01}}
\left (\frac{1}{2} \right )^{x_{\triangle}}
{ \xplus \choose \xtr}
\frac{B(\atr+\xtr,\aplus-\atr+\xplus-\xtr)}{B(\atr,\aplus-\atr)}.\nonumber
\end{eqnarray}
\end{prop}
\textit{Proof}. From (\ref{eq:pIasMixture}) we know that
$
p_H^{I}(\eta,\theta|H_0) = \sum_x p_H(\eta,\theta|x)m_{H_0}(x)
$.
Recalling (\ref{eq:quadrinomial})
and the distributional results of Lemma~\ref{lem:prior},
the expression for $p_H(\eta,\theta|x)$ is easily seen to be
$
p_H(\eta,\theta|x)
=
Beta(\eta|\atr+\xtr,\aplus-\atr+\xplus-\xtr)
\times
Beta(\theta|a_{01}+x_{01},\atr-a_{01}+\xtr-x_{01})
$
so that, conditionally on $x$, $\eta$ and $\theta$ are independent.
On the other hand, the weight $m_{H_0}(x)$ is given by
\begin{eqnarray}
m_{H_0}(x)
& = &
\int g(\xtr|\eta)h(x_{01}|\xtr,\theta=1/2)p_{H_0}(\eta) d\eta\nonumber\\
& = &
{ x_{\triangle}\choose x_{01} }
\left (\frac{1}{2} \right )^{x_{\triangle}}
\int
{\xplus \choose \xtr }
\eta^{\xtr}(1-\eta)^{\xplus-\xtr} Beta(\eta|\atr,\aplus-\atr) d\eta\nonumber \\
& = &
{x_{\triangle} \choose x_{01}}
\left (\frac{1}{2} \right )^{x_{\triangle}}
{ \xplus \choose \xtr }
\frac{B(\atr+\xtr,\aplus-\atr+\xplus-\xtr)}{B(\atr,\aplus-\atr)}.\nonumber
\end{eqnarray}
and the thesis follows.
\cvd
\smallskip

Notice that the prior $p_{H_0}(\eta)$ is set equal to $p_H(\eta)$,
i.e.~the marginal prior for $\eta$ under $H$, in~agreement with the choice that was
made in Subsection \ref{subs:single-intrinsic}. In this way we shall be able to make
more sensible comparisons with the CI-procedure.

\begin{figure}
\begin{center}
\includegraphics[width=5in]{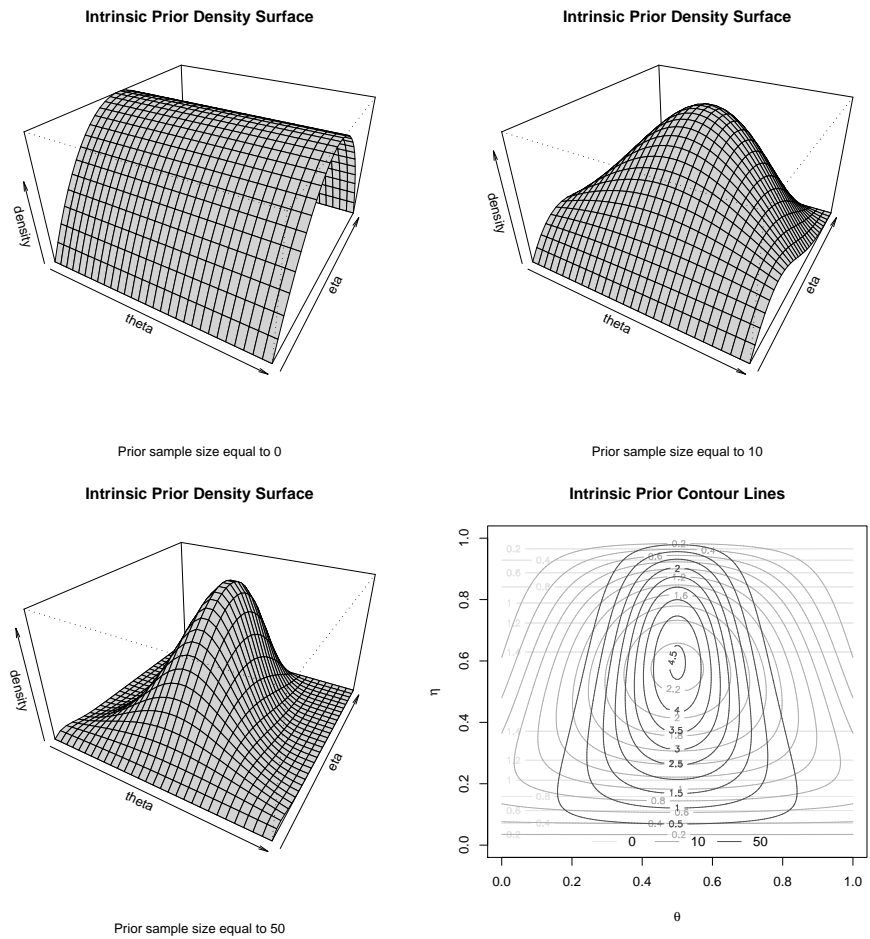}
\end{center}
\caption{Default prior density (top-left),
two different intrinsic prior densities (top-right and bottom-left),
and their contour lines (bottom-right).}
\label{fig:I prior eta theta}
\end{figure}

\begin{figure}
\begin{center}
\includegraphics[width=3.3in]{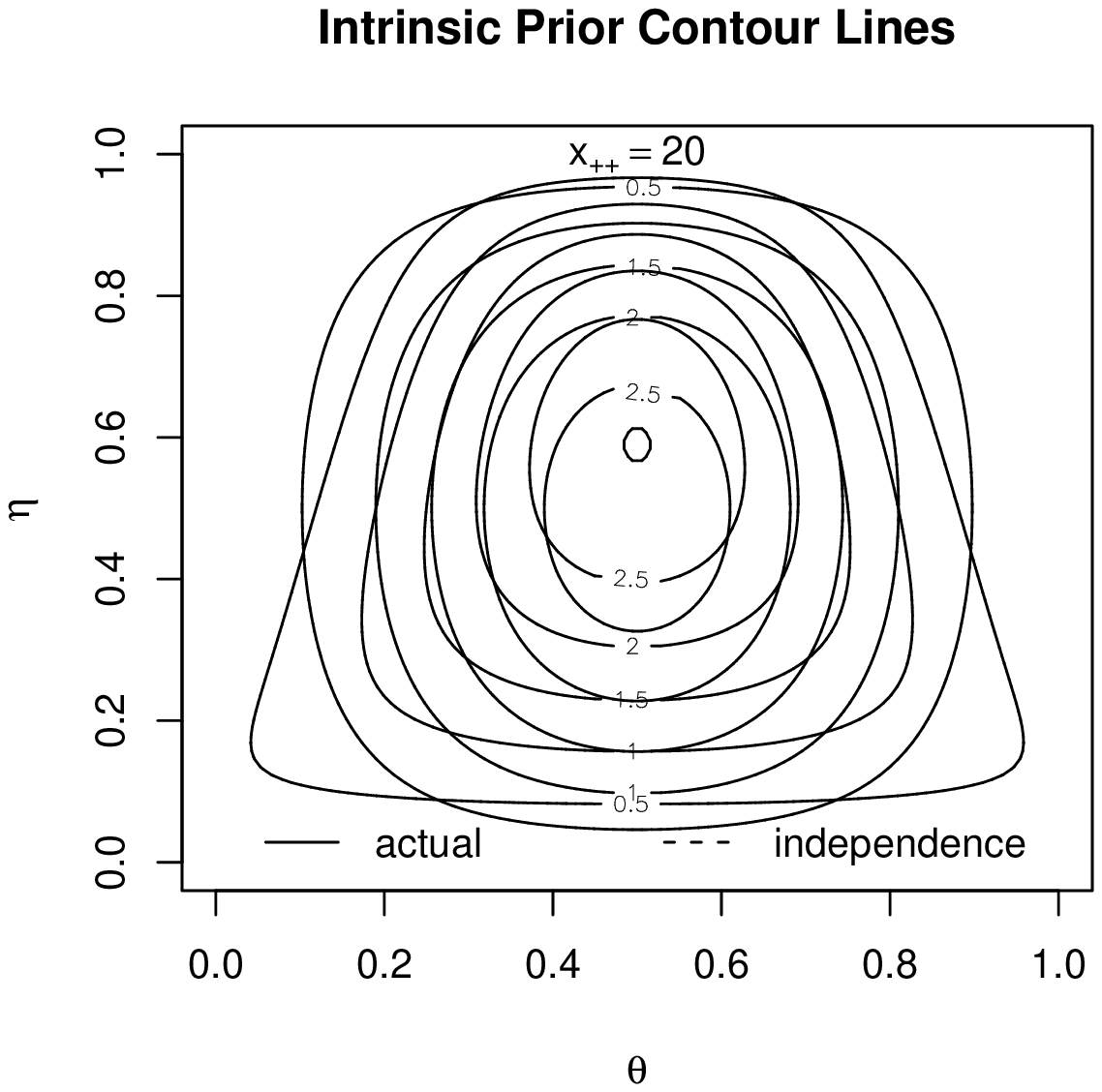}
\end{center}
\caption{Contour lines for the intrinsic prior density with $x_{++}=20$,
compared to those of the bivariate density having the same marginals with
$\eta$ and $\theta$ independent.}
\label{fig:contour lines I prior}
\end{figure}

Figure~\ref{fig:I prior eta theta} presents the default distribution
of $(\eta,\theta)$ corresponding to the uniform distribution of~$\pi$,
together with two I-priors having different prior sample sizes,
namely $\xplus=10$ and $\xplus=50$, to show the increasing concentration
of the joint distribution around the $\theta=1/2$ line.
An interesting feature is represented by the ``pear''-shaped contour lines
depicted in the bottom-right panel: they are all symmetric around the line $\theta=1/2$,
and they pile up around it as the prior sample size increases,
but the piling is not uniform with respect to~$\eta$.
More precisely, the range of credible values for $\theta$ (effective support)
becomes smaller as the value of $\eta$ increases.
This shows a dependence structure of the two parameters under the I-prior,
in contrast with the independence between $\eta$ and $\theta$ under the default prior.
Figure \ref{fig:contour lines I prior} makes this point explicit
by plotting the contour lines of the joint distribution of $(\eta,\theta)$
under the I-prior, based on $\xplus=20$ prior observations,
and those (represented by ellipsoids) of a prior having the same
marginal distributions, but embodying independence of the two parameters.

Recall that in the CI-approach parameter $\eta$ was integrated out
at the very beginning, using the same prior under $H$ and $H_0$,
since $\eta$ is a nuisance parameter. One could argue that the intrinsic approach
should implicitly recognize that $\eta$ is a nuisance parameter.
In particular, centering  the marginal prior for $\eta$ around $H_0$ makes no sense,
so that the marginal prior for $\eta$ should remain unaffected by
the intrinsic procedure. This conjecture turns out to be true,
as the following proposition shows.

\begin{prop}
The marginal distribution of $\eta$ under the intrinsic prior  is
the same as the original marginal prior, i.e.~$
p_H^{I}(\eta|H_0)=p_H(\eta)
$.
\end{prop}
\textit{Proof}.
From the expression of $p_H^{I}(\eta|H_0)$ in Proposition~\ref{prop:Iprior} we derive
\begin{eqnarray*}
p_H^{I}(\eta|H_0)
& = & \sum_{\xoi=0}^{\xplus}\sum_{\xio=0}^{\xplus-\xoi}
Beta(\eta|\atr+\xtr,\aplus-\atr+\xplus-\xtr) \\
 & \times & \frac{\xplus!}{\xtr!(\xplus-\xtr)!}
\frac{\xtr!}{\xoi!\xtr-\xoi!}
\frac{1}{2^{\xtr}}\\
& \times & \frac{B(\atr+\xoi+\xio,\aplus-\atr+\xplus-\xoi-\xio)}{B(\atr,\aplus-\atr)}\nonumber \\
& = & \sum_{\xoi=0}^{\xplus}\sum_{\xio=0}^{\xplus-\xoi}
{ \xplus  \choose \xoi \;  \xio \; (\xplus-\xtr) }
\left(\frac{\eta}{2}\right)^{\xoi}\left(\frac{\eta}{2}\right)^{\xio}(1-\eta)^{\xplus-\xtr} \\
& \times & Beta(\eta|\atr,\aplus-\atr)\nonumber \\
& = & \left[\frac{\eta}{2}+\frac{\eta}{2}+(1-\eta)\right]^{\xplus} \times p_H(\eta)
\quad=\quad p_H(\eta).
\nonumber
\end{eqnarray*}
\cvdbis
\vspace{-1.5ex}

\begin{figure}
\begin{center}
\includegraphics[width=3.3in]{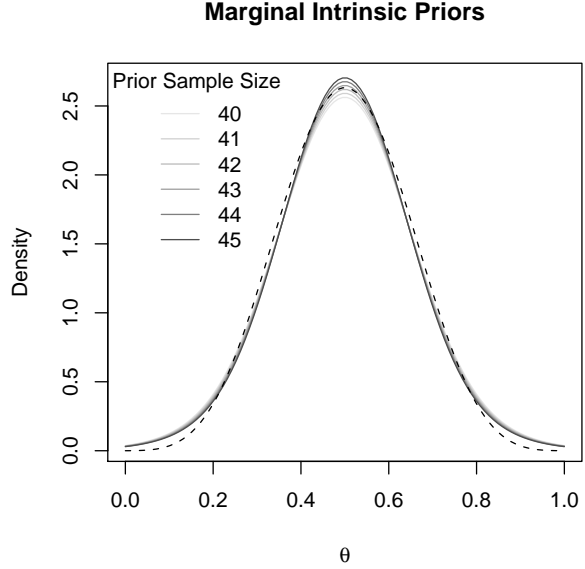}
\end{center}
\caption{Marginal intrinsic prior densities for~$\theta$
best approximating the CI-prior density based on $\xtr=20$
(dashed line).}
\label{fig:marginal and conditional I prior}
\end{figure}

On the other hand, the CI-prior for $\theta$ cannot be recovered (exactly)
within the I-procedure, as shown in Figure~\ref{fig:marginal and conditional I prior}
for a selection of prior sample sizes, $x_{++}$, best approximating the CI-prior
based on $\xtr=20$.
The marginal distribution of parameter~$\theta$ under the I-procedure is different
from the corresponding prior under the CI-procedure.
The reason lies in the different structure of the imaginary data $x$:
under the CI-procedure we have $x=(\xoi,\xtr)$,
while under the I-procedure we have  $x=(\xoi,\xtr,\xplus)$.
This gives rise to different weights:
compare $m_H(x)$ in the proof of Proposition~\ref{prop:Iprior}
to $m_H^{Co}(x)$ in the proof of Proposition~\ref{prop:CIpriorForTheta}.

Let us now consider the expression for the I-prior BF.
From Proposition \ref{prop:BFIasMixture}
we immediately derive that this is a mixture
of conditional BFs, the weights of the mixture being
the same as in the proof of Proposition~\ref{prop:Iprior}.
Hence,
$
BF^{I}_{H,H_0}(n)=\sum_x BF_{H,H_0}(n|x) m_{H_0}(x),\nonumber
$
where we need to compute
\begin{eqnarray}
BF_{H,H_0}(n|x)=\frac{\int \int g(\ntr|\eta)h(\noi|\ntr,\theta)
p_H(\eta,\theta|x)d\eta d\theta}{m_{H_0}(n)}.
\label{eq:IpriorBFconditional}
\end{eqnarray}
The denominator of~(\ref{eq:IpriorBFconditional}) is
the same as $m_{H_0}(x)$ with $x$ replaced by $n$.
On the other hand, the numerator of (\ref{eq:IpriorBFconditional})
is equal to
\begin{eqnarray}
& &
{\nplus \choose \ntr}
\frac{B(\atr+\xtr+\ntr,\aplus-\atr+\xplus-\xtr+\nplus-\ntr)}
{B(\atr+\xtr,\aplus-\atr+\xplus-\xtr)} \nonumber \\
& \times &
{\ntr \choose \noi}
\frac{B(\aoi+\xoi+\noi,\atr-\aoi+\xtr-\xoi+\ntr-\noi)}
{B(\aoi+\xoi,\atr-\aoi+\xtr-\xoi)},\nonumber
\end{eqnarray}
so that we obtain
\begin{eqnarray}
BF_{H,H_0}(n|x) m_{H_0}(x)
& = &
2^{\ntr}\frac{B(\aoi+\xoi+\noi,\atr-\aoi+\xtr-\xoi+\ntr-\noi)}
{B(\aoi+\xoi,\atr-\aoi+\xtr-\xoi)}
\nonumber\\
& \times & \frac{B(\atr+\xtr+\ntr,\aplus-\atr+\xplus-\xtr+\nplus-\ntr)}
{B(\atr+\ntr,\aplus-\atr+\nplus-\ntr)} \nonumber \\
& \times & \frac{1}{2^{\xtr}}
{\xplus \choose \xoi\; \xio\; (\xplus-\xtr)}.\nonumber
\end{eqnarray}
One can recognize that the first factor above is exactly $BF^{Co}_{H,H_0}(n|x)$,
as given in Proposition~\ref{prop:CIpriorBF}. Therefore, multiplying and dividing
by $2^{\xplus}$, we get
\begin{eqnarray}
BF^{I}_{H,H_0}(n)& = &\sum_{\xoi} \sum_{\xio}
{\xplus \choose \xoi  \xio  (\xplus-\xoi-\xio)}
\frac{1}{2^{\xtr+\xplus}} \nonumber \\
& \times &
2^{\xplus}\frac{B(\atr+\xtr+\ntr,\aplus-\atr+\xplus-\xtr+\nplus-\ntr)}
{B(\atr+\ntr,\aplus-\atr+\nplus-\ntr)} \nonumber \\
& \times & {\mbox BF^{Co}_{H,H_0}(n|x)} \nonumber.
\end{eqnarray}
Curiously, the second line in the above expression
can be further recognized as a particular BF of the hypothesis $\eta \neq \frac{1}{2}$
\textit{versus} the alternative $\eta=\frac{1}{2}$,
relative to the marginal binomial sampling model $g$,
having observed the imaginary data $x$,
and assuming a $\mbox{Beta}(\atr+\ntr,\aplus-\atr+\nplus-\ntr)$ ``prior'' distribution
on $\eta$. By writing this Bayes factor as
${\mbox BF^{M}_{\eta \neq \frac{1}{2},\eta=\frac{1}{2}}(x|n)}$,
and suitably re-writing the factor $\frac{1}{2^{\xtr+\xplus}}$,
we get to the following proposition.

\begin{prop}
\label{prop:IpriorBF}
The I-prior BF in favor of~$H$ is given by
\begin{eqnarray}
BF^{I}_{H,H_0}(n)& = &\sum_{\xoi} \sum_{\xio}
{ \xplus \choose   \xoi \;  \xio\;  \xplus-\xtr}
\left(\frac{1}{4}\right)^{\xoi} \left(\frac{1}{4}\right)^{\xio}
\left(\frac{1}{2}\right)^{\xplus-\xtr}\nonumber \\
& \times &
{\mbox BF^{M}_{ \{\eta \neq \frac{1}{2},\eta=\frac{1}{2}\}}(x|n)}
 \times  {\mbox BF^{Co}_{H,H_0}(n|x)},\nonumber
\end{eqnarray}
where
$
BF^{Co}_{H,H_0}(n|x)
=
2^{n_{\triangle}}\frac{
B(a_{01}+x_{01}+n_{01},\atr-a_{01}+x_{\triangle}-x_{01}+n_{\triangle}-n_{01})}{B(a_{01}+x_{01},\atr-a_{01}+x_{\triangle}-x_{01})}
$
and
$
\mbox BF^{M}_{\eta \neq \frac{1}{2},\eta=\frac{1}{2}}(x|n)
=
2^{\xplus}\frac{B(\atr+\xtr+\ntr,\aplus-\atr+\xplus-\xtr+\nplus-\ntr)}
{B(\atr+\ntr,\aplus-\atr+\nplus-\ntr)}.
$
\end{prop}

Hence, $BF^{I}_{H,H_0}(n)$ is a mixture of products of two BFs:
the CI-prior BF derived in the previous section,
and the marginal BF of $\eta \neq \frac{1}{2}$ \textit{versus} $\eta=\frac{1}{2}$
using imaginary data and an actual posterior distribution on $\eta$.


\section{Examples}
\label{sec:examples}

We illustrate our methodology using three examples.
They were chosen by \cite{Iron:Pere:Tiwa:2000} in order to reflect
agreement or disagreement on testing $H_0$ \textit{versus} $H$
between a frequentist approach (p-value) and a Bayesian approach
(default BF of Lemma~\ref{lem:defaultBF} derived from a uniform prior
on the vector of cell-probabilities~$\pi$).
The analysis carried out in \cite{Iron:Pere:Tiwa:2000}
also involves estimating $\theta$, as well as other parameters,
but we shall not report on it.
Conventionally, \cite{Iron:Pere:Tiwa:2000} reject the null
when the p-value is below 5\%, in the frequentist mode,
and when the default BF \textit{against} the null is above 1,
in the Bayesian mode.

We shall provide an alternative Bayesian analysis based on
our intrinsic prior methodology, also starting from
a uniform prior on the cell probabilities.
For ease of communication, we transform the BF to the probability scale of~$H_0$,
assuming prior odds equal to one,
i.e.~we consider $\prob(H_0|n)=\frac{1}{1+BF_{H,H_0}(n)}$,
and similarly for $\prob^{CI}(H_0|n)$ and $\prob^I(H_0|n)$.

For each example, rather than selecting a single intrinsic prior,
we shall assess the sensitivity of our results to prior sample size.
In the spirit of \cite{Fan:Berg:2000},
this will be achieved through a plot of the posterior probability of~$H_0$
as a function of the ratio between prior and actual sample size,
which we label~$q$. Specifically, $q=\xplus/\nplus$ for
the I-case, and $q=\xtr/\ntr$ for the CI-case, $0<q<1$.

We now present and discuss the three examples.
Although the data are fictitious, for the sake of exposition
they are presented as a survey of individuals expressing support
(``Yes'', ``No'') for the President, before and after a Presidential address.

\begin{example}
The data for this case represent a random sample of 100 interviewed individuals,
and are reported in Table \ref{tab:ex1}.
The p-value is $25\%$ and $\prob(H_0|n)=0.64$,
so that neither approach would reject $H_0$.

\begin{table}[p]
\small
\label{tab:ex1}
\begin{center}
\begin{tabular}{rrrr}
\multicolumn{4}{c}{Support for the President}\\
\hline
&\multicolumn{3}{c}{After} \\
\cline{2-4}
Before &  No & Yes & Total \\
\cline{1-1}
No     &  20 &  17 &    37 \\
Yes    &  10 &  53 &    63 \\
Total  &  30 &  70 &   100 \\
\hline
\end{tabular}
\end{center}
\caption{Data for Example~1.}
\end{table}
\end{example}

\begin{example}
In this example a random sample of 100 interviewed individuals
gave rise to Table~\ref{tab:ex2}.
The p-value is $4\%$, so that $H_0$ would be rejected,
while $\prob(H_0|n)=0.29$, leading to the same conclusion.

\begin{table}[p]
\small
\label{tab:ex2}
\begin{center}
\begin{tabular}{rrrr}
\multicolumn{4}{c}{Support for the President}\\
\hline
&\multicolumn{3}{c}{After} \\
\cline{2-4}
Before &  No & Yes & Total \\
\cline{1-1}
No     &  20 &  21 &    41 \\
Yes    &   9 &  50 &    59 \\
Total  &  29 &  71 &   100 \\
\hline
\end{tabular}
\end{center}
\caption{Data for Example~2.}
\end{table}
\end{example}

\begin{example}

In this last example the sample size is only 14
and the outcome of the survey is given by Table \ref{tab:ex3}.
The p-value is $7\%$, suggesting that $H_0$ should not be rejected;
on the other hand $\prob(H_0|n)=0.22$, pointing to the opposite conclusion.
Thus, for this example, there seems to be a clear conflict
between the conclusions under the two approaches.

\begin{table}[p]
\small
\label{tab:ex3}
\begin{center}
\begin{tabular}{rrrr}
\multicolumn{4}{c}{Support for the President}\\
\hline
&\multicolumn{3}{c}{After} \\
\cline{2-4}
Before &  No & Yes & Total \\
\cline{1-1}
No     &   1 &   7 &    8 \\
Yes    &   1 &   5 &    6 \\
Total  &   2 &  12 &   14 \\
\hline
\end{tabular}
\end{center}
\caption{Data for Example~3.}
\end{table}

\end{example}

\begin{figure}
\label{intriPriorAndPost}
\begin{center}
\includegraphics[width=5in]{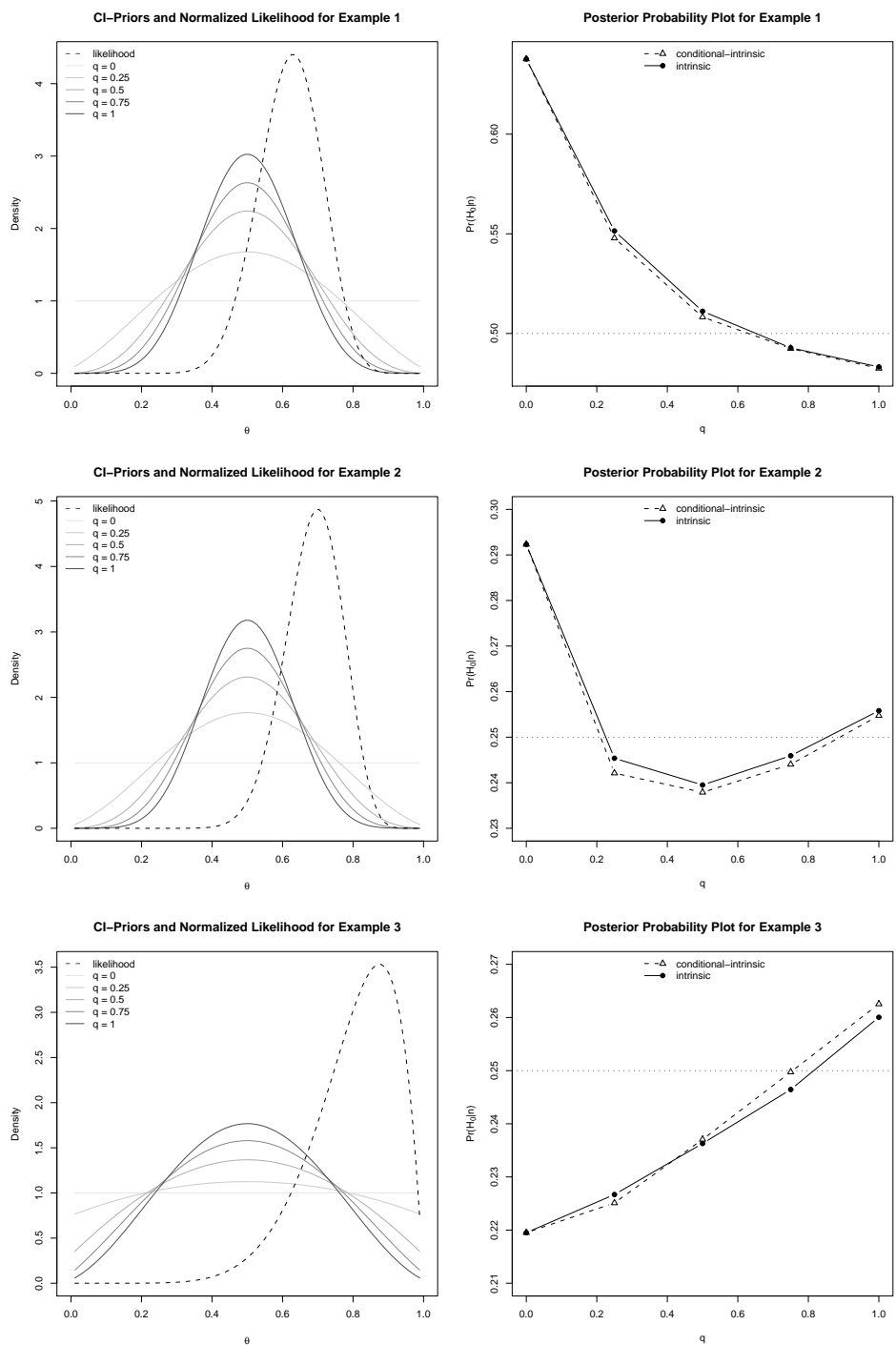}
\end{center}
\caption{Intrinsic prior and posterior probability of~$H_0$:
sensitivity to prior sample size for each of the three examples
($q$ is the fraction of prior to actual sample size).}
\end{figure}

We now proceed with the intrinsic analysis of the three examples,
having in  mind that sensitivity to prior specifications is a major concern when evaluating BFs. Accordingly, our results  are visually summarized
in Figure~\ref{intriPriorAndPost}. For each example we present two panels:
on the left-hand-side we plot the CI-prior density
for selected values of the fraction of prior to actual sample size $q$,
namely $q=0,0.25,0.50,0.75,1$.
Clearly, $q=0$ corresponds to the starting
default prior for $\theta$, that is the uniform distribution
on the interval $(0,1)$. We also plot in this panel
the normalized likelihood function.
In the right panel we present the plots of
$\prob^{CI}(H_0|n)$ and $\prob^I(H_0|n)$,
as a function of~$q$.
Notice that for all three datasets
the two plots are remarkably similar.

Consider first Example 1. The left panel reveals a mild conflict between
the likelihood and the null hypothesis, as well as the CI-prior
which by construction is centered on $H_0$. The right panel shows
that $\prob^{CI}(H_0|n)$ and $\prob^I(H_0|n)$ are decreasing
for $0<q<1$, spanning
a range from 0.64 (default value) to~0.47.
This can be explained as follows: the data do not contradict explicitly $H_0$
and thus there is scope for diminishing the probability of~$H_0$
by peaking the prior under~$H$ around~$H_0$.
Although the curves cross the 50\% threshold, this only occurs
if $q$ exceeds 70\%: a~value that appears to overemphasize the role of the prior.
In conclusion, the rejection of~$H_0$ on conventional grounds
is not totally robust, but is broadly justified.

Consider now Example 2. The left panel reveals an appreciable conflict
between the likelihood on the one hand, and the null and the CI-prior on the other.
This aspect translates into a posterior probability of $H_0$
approximately ranging from~0.24 to~0.29 (default value).
Notice that in this example the behavior of $\prob^{CI}(H_0|n)$ and $\prob^I(H_0|n)$
is not monotone as a function of~$q$. To understand why this occurs,
recall that the data are in disagreement with the null.
As a consequence, centering the prior under~$H$ around the null
can only benefit the posterior probability of~$H_0$ provided
the prior sample size is not too large;
otherwise the prior under~$H$ will be too similar to that under~$H_0$.
%
In conclusion, the conventional rejection of~$H_0$ appears fully robust;
yet the default value of $\prob(H_0|n)$ is possibly too high,
since a moderate prior sample size determines a relatively sharp decrease
in $\prob^{CI}(H_0|n)$ and $\prob^I(H_0|n)$. A value close to
the 0.25 threshold seems therefore more reasonable.

Finally, turn to Example 3. Here the likelihood reveals a marked disagreement
with the null, although it has a rather heavy left tail because of the
limited number of observations.
The family of CI-priors is accordingly more dispersed around the value $\theta=1/2$,
because they are based on prior sample sizes which are a fraction of $\ntr=8$.
Notice that $\prob^{CI}(H_0|n)$ and $\prob^I(H_0|n)$ are now monotone increasing
from~0.22 (default value) to~0.26.
This is a clear example where the intrinsic prior methodology contributes to
water down the evidence against the null, thus defying its original purpose.
This happens because the data are already in good agreement with the alternative,
and thus peaking the prior under~$H$ around~$H_0$ only accumulates probability mass
in areas not supported by the likelihood, thus increasing
the posterior probability of $H_0$.
Hence, in this case too, the conventional rejection of~$H_0$ is robust.

\section{Discussion}

In this paper we considered Bayesian hypothesis testing
for the comparison of two nested hypotheses, using the BF as a measure of evidence.
Alternative approaches may be useful,
e.g.~those based on measures of divergence such as the symmetrized Kullback-Leibler;
see \citet{Bern:2005}. However, when converted to the probability scale of the null,
assuming prior odds equal to one, the BF is a communication device
of unsurpassed clarity.
However, since the BF is strongly influenced by  prior specification,
whose effect does not decrease as the sample size increases,
sensitivity analysis is needed to evalutate robustness of the conclusions.

An objective framework for sensitivity analysis is offered
by the intrinsic prior methodology: starting from default priors
both under the null and the alternative (encompassing) hypothesis,
a family of prior distributions under the alternative is constructed
such that its elements are increasingly concentrated around the null,
as prior sample size grows. In this way, we avoid to unduly favor
the null, especially when the data do not contradict it explicitely

\begin{figure}
\begin{center}
\includegraphics[width=3.3in]{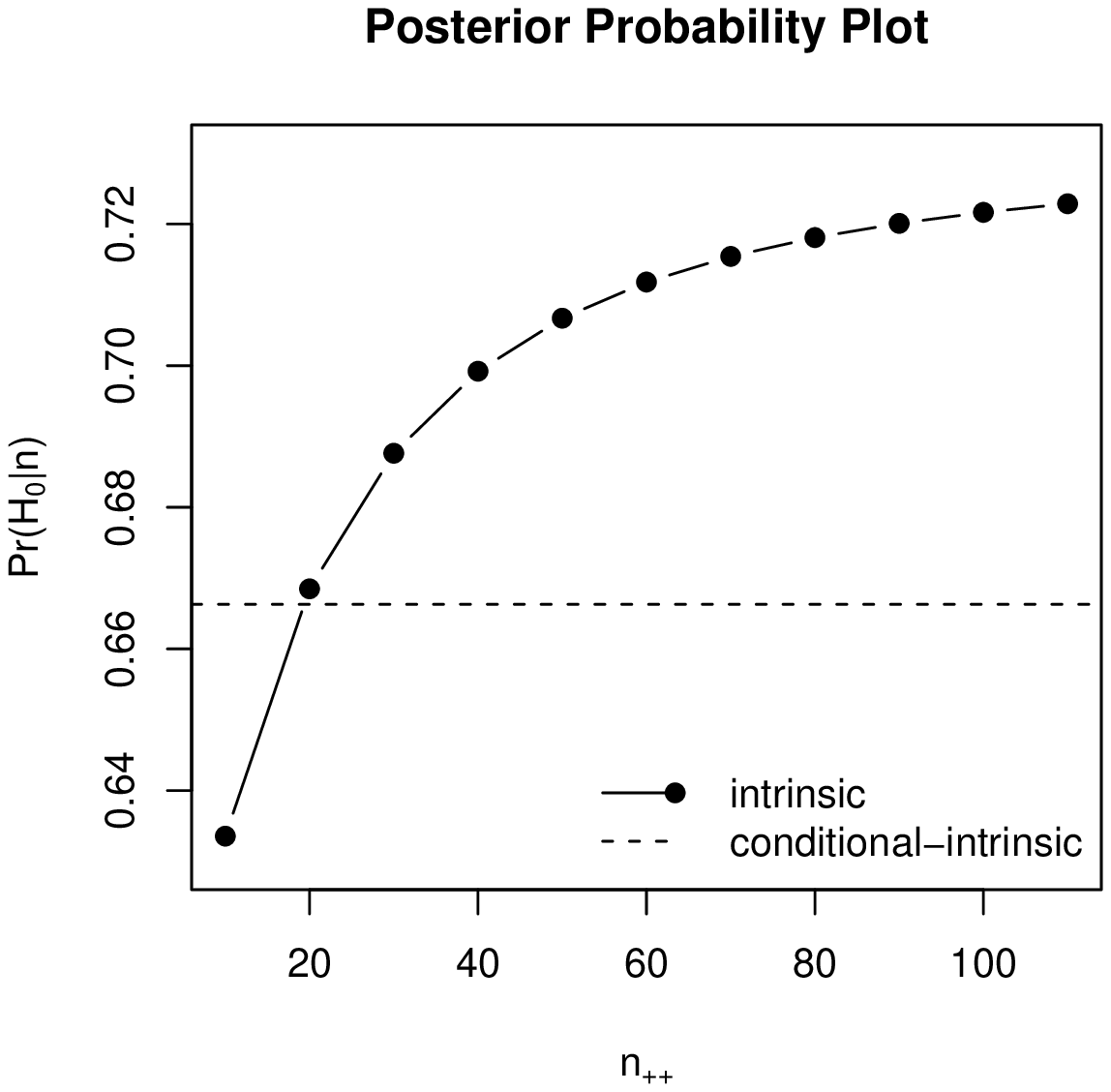}
\end{center}
\caption{Posterior probability of the null, as a function of~$\nplus$,
for a table~$n$ such that $\noi=\nio=5$, with $\xplus=10$ and $\xtr=5$ for I-priors
and CI-priors, respectively.}
\label{fig:dependOnNPP}
\end{figure}

For the problem of testing the equality of two correlated proportions,
there is a choice between using the full sampling model~(\ref{eq:quadrinomial})
and the conditional sampling model given the overall number of swing.
We considered both approaches, obtaining I-priors and CI-priors, respectively.
When applied to the examples of Section~\ref{sec:examples},
the two approaches proved to be essentially equivalent,
altough in general they are not. In particular, as shown in Figure~\ref{fig:dependOnNPP},
the posterior probability of the null depends on the total sample size~$\nplus$
in the I-approach only; see \citet{Ghos:Chen:Ghos:Agre:2000} for similar results
obtained using a Rash-type Bayesian hierarchical model.

A point that was not addressed in this paper concerns consistency of
the CI-prior and I-prior Bayes factors. \cite{Case:More:2006},
discussing tests of independence in two-way contingency tables,
show that the intrinsic BF is consistent provided the ratio of prior
to (actual) sample size goes to zero as the  sample size goes to infinity.

\pagebreak

Finally, the problem of comparing two correlated proportions
is a special case of the more general problem of testing the hypothesis
of marginal homogeneity in square contingency tables;
see \citet[Ch.~3]{Berg:1997} for an extensive review.

\section*{Acknowledgements}

Research partially supported by MIUR, Rome (PRIN 2005132307)
and  the University of Pavia. The seminal idea of this paper was developed
while the Authors were visiting CRiSM and the Department of Statistics
at the University of Warwick, England. We thank the above institutions
for hospitality and support. The authors would like to thank
Jim Berger, El\'\i as Moreno and Mario Peruggia for helpful comments.

\bibliographystyle{biometrika}
\bibliography{marmod}

\end{document}